\documentclass[a4paper,10pt]{article}
\usepackage{paper}

\def\thetitle{About every convex set in any\\ generic Riemannian manifold}
\hypersetup{
pdftitle={\thetitle},
pdfauthor={Alexander Lytchak and Anton Petrunin}
}

\begin{document}

\title{\thetitle}
\author{Alexander Lytchak and Anton Petrunin}
\date{}
\maketitle

\begin{abstract}
We give a necessary condition on a geodesic in a Riemannian manifold that can run in some convex hypersurface.
As a corollary, we obtain peculiar properties that hold true for \emph{every} convex set in any \emph{generic} Riemannian manifold $(M,g)$.
For example, if a convex set in $(M,g)$ is bounded by a smooth hypersurface, then it is strictly convex.
\end{abstract}







\section{Introduction}
Let $\mathfrak{C}$ be the convex hull of a compact subset $Q$ in the Euclidean space $\R^m$.
By Carathéodory's theorem \cite{Handbook}, $\mathfrak{C}$ is the set of all convex combinations of at most $(m+1)$-tuples of points on $Q$.
Thus, $\mathfrak{C}$ is a compact convex subset.  Any point in $\mathfrak{C}\setminus Q$ is an inner point of a line segment contained in $\mathfrak{C}$;
that is, the complement $\mathfrak{C} \setminus Q$ does not contain extreme points of $\mathfrak{C}$. 

The compactness of the convex hull and, therefore, the existence of a huge variety of convex subsets with many non-extreme points on the boundary, admits a straightforward generalization to the sphere and the  Lobachevsky space; moreover, it 
holds \emph{locally} in any two-dimensional Riemannian manifold.

Recall that a set $\mathfrak{C}$ in a Riemannian manifold $(M,g)$ is called \emph{convex} if for any pair of points $x,y\in \mathfrak{C}$ any minimizing geodesic $[x,y]$ lies in $\mathfrak{C}$.
A point in $\mathfrak{C}$ is called \emph{extreme} if it does not lie in an interior of a geodesic in~$\mathfrak{C}$.
The \emph{convex hull} of a set $Q\subset M$ is the minimal convex subset of $M$  that contains $Q$.

It seems to be a folklore belief that a version of the statement above should hold true in all Riemannian manifolds;
see the discussion at mathoverflow \cite{petrunin-2009}.
In the present note we prove that the somewhat counter-intuitive opposite is the case for \emph{generic} Riemannian manifolds.
It agrees with the pattern: \emph{a typical object in your favorite theory looks like nothing you have ever seen before}.

Further Riemannian manifolds will be assumed to be connected and $\mathcal C^\infty$-smooth.
Given a positive integer $k$, we say that a property $\mathcal P$ holds for \emph{$\mathcal C^k$-generic} Riemannian metric $g$ on a manifold $M$ 
if the property $\mathcal P$ holds for a dense \emph{G-delta set} (that is, a countable intersection of open subsets) of metric tensors in the $\mathcal C^k$-topology.

\begin{thm}{Main theorem}\label{thm:main}
Let $\mathfrak C$ be an arbitrary convex subset of a $\mathcal C^2$-generic Riemannian manifold $(M,g)$.
Then the set of non-extreme points in $\mathfrak C$ is the union of an open set and at most countable family of geodesics in $(M,g)$. 

In particular, if $\dim M\ge 3$ and no connected component of $\mathfrak C$ is  a geodesic, then the set of
extreme points of $\mathfrak{C}$ is dense in~$\partial\mathfrak{C}$.
\end{thm}

Note that our definition of convexity does not require connectedness. However, any convex subset $\mathfrak C$ is locally connected, and $\mathfrak C$  is 
connected if the manifold $M$ is complete or $\mathfrak C$ is contained  in some compact convex subset of $M$.

If $\dim M =2$, the statement is rather trivial and holds true for \emph{all} Riemannian metrics not only the generic ones.
As a consequence  of the main theorem for $\dim M \geq 3$, we obtain the following:

\begin{thm}{Corollary}\label{cor:caratheodory}
Let $Q$ be a closed subset of a $\mathcal C^2$-generic Riemannian manifold $(M,g)$ of dimension at least $3$. If $Q$ does not lie in a  geodesic and
the convex hull $\mathfrak C$ of $Q$ is closed and connected, then $\partial \mathfrak C \subset Q$.
\end{thm}

The corollary gives a positive resolution of a conjecture formulated by Marcel Berger
\cite[Note 6.1.3.1]{berger-2003}, stating that convex hulls of $3$ points in most Riemannian manifolds do not need to be  compact.
Probably  the following more exact form of Berger's conjecture might be squeezed out from our key lemma.%
\footnote{More open questions are listed in  Appendix \ref{app:remarks}.}

\begin{thm}{Conjecture}  Let $(M,g)$ be an arbitrary Riemannian manifold of dimension at least 3.
If the convex hull of any 3-point subset is compact,
then $(M,g)$ has constant curvature.
\end{thm}

The following corollary is essentially known, for $\dim M=3$ its proof has been sketched by Robert Bryant \cite{Bryant} and, for $\dim M\geq 4$, it was proved by Thomas Murphy and Frederick Wilhelm \cite{Wilhelm}.

\begin{thm}{Corollary}\label{cor:main}
Let $(M,g)$ be a $\mathcal C^2$-generic Riemannian manifold.
Then any connected  convex subset $\mathfrak C$ of $(M,g)$ is either contained in a geodesic
or \emph{full-dimensional}; that is, the interior of $\mathfrak C$ is nonempty.
\end{thm}

The proofs are built on the following proposition.
Its formulation uses the notion of \emph{rank} of a point $p$ in a closed convex set $\mathfrak{C}$;
we define it as the dimension of the maximal linear subspace in the tangent cone to $\mathfrak{C}$ at $p$.

\begin{thm}{Main proposition}\label{prom:rank}
Suppose $\mathfrak{C}$ is a closed convex set in a $\mathcal C^2$-generic $m$-dimensional Riemannian manifold $(M,g)$.
Then all non-extreme points of $\mathfrak{C}$ have rank either $1$ or $m$.

In particular, if $\dim M\ge 3$ and $\mathfrak{C}$ is bounded by a $\mathcal{C}^1$-smooth hypersurface, then $\mathfrak{C}$ is \emph{strictly convex};
that is, all boundary points of $\mathfrak{C}$ are extreme.
\end{thm}

The proof relies on the key lemma stated in the following section;   
it describes a necessary condition on a geodesic in a Riemannian manifold that stays in convex subset $\mathfrak C$.
If the geodesic  lies in $\partial \mathfrak C$ and  contains a point of rank at least 2, then this condition implies a non-trivial
property  of the curvature tensor.
Then we show that the curvature tensor of a generic Riemannian manifold does not meet this property.
The latter part is technical but straightforward; it is done by applying the Thom transversality theorem; see Appendix \ref{sec:normalization}.

\parbf{Acknowledgments.}
We thank Mohammad Ghomi and Frederick Wilhelm for their interest in our result,
the anonymous referee for helpful criticism.
Alexander Lytchak was partially supported by the DFG grant, no. 281071066, TRR 191.
Anton Petrunin was partially supported by the NSF grant, DMS-2005279.

\section{Key lemma}\label{sec:key}

Let $\mathfrak{C}$ be a closed convex set in an $m$-dimensional Riemannian manifold $(M,g)$.
Recall that $\T_x=\T_xM$ denotes the \emph{tangent space} of $M$ at $x$.
The \emph{tangent cone} $\K_x=\K_x\mathfrak{C}\subset \T_x$ at $x\in\mathfrak{C}$ is defined as the closure of the set of all velocity vectors of geodesics that start at $x$ and run in $\mathfrak{C}$.

Given $x\in \mathfrak{C}$, denote by $\L_x=\L_x\mathfrak{C}$ the \emph{maximal linear subspace} of $\K_x$.
We define the \emph{rank} of $x$ in $\mathfrak{C}$ as the dimension of $\L_x$.

Note that $\K_x$ is a convex cone in $\T_x$; in particular, $\L_x=\K_x\cap (-\K_x)$.
Further $\K_x$ coincides with 
$\T_x$ if and only if
a neighborhood of $x$ lies in the interior of $\mathfrak{C}$.
In other words, $x$ has rank $m$ if and only if $\mathfrak{C}$ contains
a neighborhood of $x$.

Given a tangent vector $\vec x\in\T_pM$, consider the  \emph{Jacobi operators} of order $k$
\[R^k_\vec x\:\vec{v}\mapsto \nabla^{k-2}_\vec x\Rm(\vec{v},\vec x)\vec x,\]
where $\Rm$ denotes the curvature tensor of $g$;
we set $R^1=0$.
Note that (i) $R^k_\vec x\:\T_p\z\to \T_p$ is a self-adjoint operator, (ii) $\vec x\mapsto R^k_\vec x$ is a homogeneous polynomial of degree $k$, and (iii) 
\[R^k_\vec x\cdot\vec x=0 \eqlbl{eq:RXX=0}\]
for any $k$ and $\vec x\in\T_p$.
The Jacobi equation along a geodesic $\gamma$ takes the form 
\[\nabla^2_{\gamma'}\cdot\vec i+R^2_{\gamma'}\cdot \vec i=0.\]

\begin{thm}{Key lemma}\label{lem:key}
Let $(M,g)$ be a Riemannian manifold and $\gamma\:(a_0,b_0)\to M$ be a geodesic that runs in a closed convex set $\mathfrak{C}\subset (M,g)$.
Then the tangent cones of $\mathfrak{C}$ are parallel along $\gamma$; that is, the parallel translation along $\gamma$ defines a bijection between the tangent cones $\K_{\gamma(a)}\mathfrak{C}$ and $\K_{\gamma(b)}\mathfrak{C}$ for any $a,b \in (a_0,b_0)$.

Moreover, for any $a\in (a_0,b_0)$ the following conditions  hold:
\begin{subthm}{lem:key:a}
For any $\vec{v}\in \K_{\gamma(a)}\mathfrak{C}$ we have
\[R^2_{\gamma'(a)}\cdot \vec{v}\in \K_\vec{v}[\K_{\gamma(a)}\mathfrak{C}].\]
\end{subthm}

\begin{subthm}{lem:key:b} 
$\L_{\gamma(a)}\mathfrak{C}$ is an invariant subspace of $R^2_{\gamma'(a)}\:\T_{\gamma(a)}\to\T_{\gamma(a)}$.
\end{subthm}

\end{thm}

The proof uses the fact that the parallel translation can be defined via geodesics.
In a similar way, this observation was used in \cite[Section 13]{Ber-Nik} and~\cite{Petruninpar}.
In fact the main part of the key lemma follows from
\cite{Petruninpar}.

\parit{Proof of \ref{lem:key}.} 
Since all statements are local, we may replace $(M,g)$ by its small open convex subset.
By doing so we may assume that any pair of points of $(M,g)$ is connected by a unique geodesic and there are no conjugate points.
In particular, for any subinterval $[a,b]\subset (a_0,b_0)$ and any tangent vectors $\vec{v} \in \T_{\gamma (a)}$ and $\vec{w} \in \T_{\gamma (b)}$ there exists a unique Jacobi field $\vec i$ along $\gamma$ 
such that $\vec i(a)=\vec{v}$ and~$\vec i(b)=\vec{w}$.

Since Jacobi fields are variational fields of geodesic variations, 
the convexity of $\mathfrak{C}$ implies the following: 

\begin{thm}{Observation}
Suppose $\vec i$ is a Jacobi field along 
$\gamma$ and $a_0<a<t<b<b_0$.
If 
$\vec i(a)\in \K_{\gamma(a)}\mathfrak{C}$ and $\vec i(b)\in \K_{\gamma(b)}\mathfrak{C}$,
then $\vec i(t)\in \K_{\gamma(t)}\mathfrak{C}$.
\end{thm}

Choose a subinterval $[a,b] \subset (a_0,b_0)$.
Given a large positive integer $k$, consider the arithmetic progression
$t_0,\dots,t_{k+1}$ such that $t_0=a$ and $t_k=b$.

Choose a tangent vector $\vec{v}_0\in\T_{\gamma(a)}$.
Consider the sequence of vectors $\vec{v}_i\z\in\T_{\gamma(t_i)}$ defined recursively by $\vec{v}_{i+1}=2\cdot \vec i_i(t_{i+1})$, where $t\mapsto \vec i_i(t)$ denotes the Jacobi field along $\gamma$ such that $\vec i_i(t_i)=\vec{v}_i$ and $\vec i_i(t_{i+2})=0$.

\begin{figure}[ht!]\vskip-0mm\centering\includegraphics{mppics/pic-1}\end{figure}

Define $\iota_k\:\T_{\gamma(a)}\to \T_{\gamma(b)}$ by setting $\iota_k(\vec v_0)\df \vec v_k$.
According to the observation, if $\vec v_0\in \K_{\gamma(a)}\mathfrak{C}$, then $\iota_k(\vec v_0)\in \K_{\gamma(b)}\mathfrak{C}$.
As observed in \cite{Ber-Nik} and~\cite{Petruninpar}, $\iota_k(\vec v_0)$ converges to the parallel translation of $\vec v_0$ along $\gamma$ as $k\to \infty$.
Since $\K_{\gamma(b)}\mathfrak{C}$ is closed,
the parallel translation along $\gamma$ maps $\K_{\gamma(a)}\mathfrak{C}$ in $\K_{\gamma(b)}\mathfrak{C}$.
Switching the direction of $\gamma$, we get the opposite inclusion.
That is, the tangent cones $\K_{\gamma(t)}\mathfrak{C}$ are parallel along $\gamma$ --- the main part is proved.

Let us use the parallel translation along $\gamma$ to identify the tangent spaces at points on $\gamma$.
This way we identify the tangent cones $\K_{\gamma(t)}\mathfrak{C}$ for all $t$;
denote the obtained cone by $\K$.

For $\vec{v}\in \K$ and small $\epsilon>0$, consider the unique Jacobi field $\vec i_\epsilon$ along $\gamma$ with $\vec i_\epsilon (a+\eps)\z=\vec i_\epsilon(a-\eps)=\vec{v}$.
Due to the Jacobi equation,
\[\vec i_\epsilon (a)=\vec{v} +\tfrac12\cdot\eps^2\cdot R^2_{\gamma'(a)}\cdot \vec{v} +o(\eps^2).\]
According to the observation, $\vec i_\epsilon(a)\in \K$ for any $\eps>0$.
Since $\K$ is a closed convex cone, we get $R^2_{\gamma'}\cdot \vec{v}\z\in \K_\vec{v}\K$ --- \ref{SHORT.lem:key:a} is proved.

Finally,
$\vec{v}\in \L_{\gamma(a)}\mathfrak{C}$ $
\iff$
$\vec{v}, -\vec{v}\in \K$
$\iff$
$\K_\vec{v}\K=\K$.
Therefore, if $\vec{v}\in \L_{\gamma(a)}\mathfrak{C}$, then $\pm R^2_{\gamma'(a)}\cdot \vec{v}\in \K$, and hence $R^2_{\gamma'(a)}\cdot \vec{v}\in \L_{\gamma(a)}\mathfrak{C}$.
That is, $\L_{\gamma(a)}\mathfrak{C}$ is an invariant subspace of $R^2_{\gamma'(a)}$ --- \ref{SHORT.lem:key:b} is proved.
\qeds

\section{Main proposition}

In this section we will prove the main proposition \ref{prom:rank} modulo one claim; let us introduce notations to state it.

Let $M$ be a smooth $m$-dimensional manifold with a Riemannian metric $g$.
Suppose $\vec x$ is a nonzero tangent vector at a point $p\in M$.
Recall that $R^k_\vec x\:\T_p\z\to\T_p$ denotes the Jacobi operators of $g$ of order $k$ for a tangent vector $\vec x\in\T_p$.
An invariant subspace $V\subset \T_p$ of $R^k_\vec x$ will be called \emph{exceptional} if $V\ni \vec x$ and $1< \dim V<m$.
(Recall that $R^k_\vec x\cdot \vec x=0$
for any $k$ and $\vec x\in \T_p$.
Therefore, the subspace spanned by $\vec x$ is always an invariant subspace of $R^k_\vec x$ for any $k$.)

We will  say that a metric $g$ on a manifold $M$ is $k$-exceptional if there exists a point $p\in M$ and a non-zero vector $\vec x\in T_p M$,
such that the operators   $R^2_\vec x,\dots, R^k _\vec x$ have a common exceptional invariant subspace.  

\begin{thm}{Claim}\label{clm:codim-sigma} 
For any smooth manifold $M$, there exists an integer $k$ such that the 
$\mathcal C^k$-generic Riemannian metric is not $k$-exceptional.
\end{thm}

For $k=2$  (and, probably,  also for $k=3$) every Riemannian metric is $k$-exceptional.
However, for larger $k$, the
$k$-exceptionality defines more and more restrictions on the curvature tensor.
Therefore, it is not surprising that \emph{most}
Riemannian metrics are not $k$-exceptional, for sufficiently large $k$.
A formal proof of this claim is built on 
Thom transversality theorem;
it will be derived in Appendix~\ref{sec:normalization}.

\parit{Proof of \ref{prom:rank} modulo \ref{clm:codim-sigma}.}
Suppose that $p$ is a nonextreme point of $\mathfrak{C}$;
that is, $p$ lies on a nonconstant geodesic $\gamma\:(a,b)\to\mathfrak{C}$.

According to the key lemma (\ref{lem:key}), the family of \emph{maximal linear subspaces} $\L_{\gamma(t)}\mathfrak{C}$ of $\K_{\gamma(t)}\mathfrak{C}$ is parallel along $\gamma$ and invariant for $R^2_{\gamma(t)}$.
Note that $\L_p$ is exceptional if and only if the rank of $p$ is neither $1$, nor $m$.

Further, if a nontrivial geodesic $\gamma$ admits a parallel family $\L_t\subset \T_{\gamma(t)}$ of exceptional invariant subspaces for all $R^2_{\gamma(t)}$, then we say that $\gamma$ is \emph{exceptional}.
So, it is sufficient to show that $\mathcal C^2$-generic Riemannian manifolds $(M,g)$ do not have exceptional geodesics.

Choose a compact subset $K\subset M$ and $\eps>0$.
Consider the set $Z(K,\eps)$ of all Riemannian metrics $g$ on $M$ such that there exists an exceptional geodesic $\gamma$ in $(M,g)$ that starts at a point in $K$ and has length $\eps$.
Observe that the geodesics and the curvature tensor depend continuously on the Riemannian metric in $\mathcal C^2$-topology.
Therefore, the set $Z(K,\eps)$ is closed with respect to the $\mathcal C^2$-topology.

Suppose $\gamma$ is an exceptional geodesic that passes thru $p$ in the direction $\vec x$.
By taking covariant derivatives along $\gamma$, we get that the Jacobi operators $R^k_\vec x$ have a common exceptional invariant subspace $\L_p$, for all $k \geq 2$.
In other words, for any integer $k \geq 2$ we have
\[Z^k(K)\supset Z(K,\eps),\eqlbl{eq:Zk<Z(K,eps)}\]
where $Z^k(K)$ denotes the set of all smooth Riemannian metrics on $M$ such that for some $p\in K$ and $\vec x\in\T_p\setminus \{0 \}$ the operators $R^2_\vec x,\dots, R^k _\vec x$ have a common exceptional invariant subspace.

By the very definition of $Z^k(K)$, it is closed with respect to $\mathcal C^{k}$-topology on the space of all Riemannian metrics on $M$.
By Claim~\ref{clm:codim-sigma}, we can choose $k$ so that  $Z^k(K)$ is $\mathcal{C}^k$-meager for any $K$;
that is, its complement is a dense G-delta set in the space of all Riemannian metrics on $M$ with $\mathcal{C}^k$-topology.

Since $Z(K,\eps)$ is closed with respect to the $\mathcal C^2$-topology, \ref{eq:Zk<Z(K,eps)} implies that $Z(K,\eps)$ is $\mathcal{C}^2$-meager in the space of all Riemannian metrics on $M$.

Choose a nested sequence of compact sets $K_1\z\subset K_2\subset \dots$ that cover $M$ and set $\eps_n=\tfrac1n$.
Set 
\[Z(M)=\bigcup_n Z(K_n,\eps_n);\]
since $Z(K_n,\eps_n)$ is $\mathcal{C}^2$-meager for every $n$, so is $Z(M)$.

It remains to note that $g\in Z(M)$ if and only if $(M,g)$ has an exceptional geodesic.
\qeds

\section{Main theorem}

The following proposition is a special case of a result of Nan Li and Aaron Naber \cite[Theorem 1.6]{li-naber}.
It also can be deduced from the result of Luděk 
Zajíček~\cite{zajicek}.

\begin{thm}{Proposition}\label{prop:rectifiable}
Let $\mathfrak{C}$ be a closed convex set in a Riemannian manifold $(M,g)$.
Then the set of points in $\mathfrak{C}$ with rank at most $k$ is countably \emph{$k$-rectifiable};
that is, this set can be  covered by images of a countable set of Lipschitz maps $\RR^k\to (M,g)$.
In particular, this set contains at most countably many disjoint Borel sets with positive $k$-dimensional Hausdorff measure.
\end{thm}

\parit{Proof of \ref{cor:main} and \ref{thm:main}.}
We may assume that $\mathfrak {C}$ is connected. Assume, in addition, that  $\mathfrak{C}$ is closed.

According to \cite[Theorem 1.6]{cheeger-gromoll}, a connected closed convex set $\mathfrak{C}$ in a Riemannian manifold $(M,g)$ is homeomorphic to a manifold with boundary, say~$\mathfrak{B}$.
Moreover, the complement $\mathfrak{C}\backslash \mathfrak{B}$ is a totally geodesic submanifold of $(M,g)$; denote its dimension by $d$.

The tangent cone $\K_p \mathfrak C$ at any $p \in \mathfrak C\setminus \mathfrak B$ is a $d$-dimensional linear space.
By the main proposition (\ref{prom:rank}), $d=0, 1$, or $m$.
If $d=0$, then $\mathfrak{C}$ is a single point.
If $d=1$, then $\mathfrak C\setminus \mathfrak B$ is a  geodesic in $(M,g)$;
hence $\mathfrak C$ is contained in a  geodesic as well.
If $d=m$, then by the invariance of domain we have  $\mathfrak C\setminus \mathfrak B$ is open in $M$; that is, $\mathfrak C$ is full-dimensional --- \ref{cor:main} is proved.

By the main proposition (\ref{prom:rank}) any non-extreme point $x\in \partial \mathfrak C$ has rank 1.
Thus, there is a unique line in $\K_x\mathfrak C$ and it is the tangent line of a geodesic $\gamma\subset\mathfrak C$ that has $p$ as an inner point. 

Let us extend $\gamma$ to a maximal open interval so that $\gamma$ stays in $\mathfrak C$;
note that $p$ uniquely defines $\gamma$. 
By the main statement of the key lemma, all points on $\gamma$ lie on $\partial \mathfrak C$.
By definition, all such geodesics consist of non-extreme points.

It gives a subdivision of non-extreme points of $\partial\mathfrak C$ into geodesics with positive lengths.
By \ref{prop:rectifiable}, there are only countably many such geodesics.

If $\mathfrak C$ is not closed, consider its closure $\bar {\mathfrak C}$;
denote by $\bar{\mathfrak{B}}$ its boundary.
Note that $\bar{\mathfrak C}$ is locally convex and the above arguments apply to closed locally convex subsets without changes.
Observe that any nonextreme point of $\mathfrak{C}$ is a  nonextreme point of $\bar{\mathfrak{C}}$ and $\bar{\mathfrak{C}}\backslash\bar{\mathfrak{B}}\subset\mathfrak{C}$ \cite[Lemma 1.5]{cheeger-gromoll}.
Hence, the statement follows.
\qeds

\parbf{Few words before the proof of  \ref{cor:caratheodory}.}
Let $Q$ be a subset of a Riemannian manifold $M$.
Set $Q=Q_0$ and let inductively $Q_{i+1}$ to be the union of all minimizing geodesics between pairs of points of $Q_i$.
By definition, the increasing countable union $\mathfrak{C}= \bigcup_iQ_i$ is the convex hull of $Q$.  
By this description, any point in $\mathfrak{C} \setminus Q$ is a non-extreme point of the convex set $\mathfrak {C}$.

Note, that if $M$ is complete and $Q$ is compact, then each $Q_i$ is  compact.
In the Euclidean space $M=\R^m$  (as well as in the round sphere or in the Lobachevsky space) Carath\'eodory's theorem \cite{Handbook} implies $\mathfrak{C} =Q_m$.
As a consequence of Corollary \ref{cor:caratheodory}, we will know that in a generic Riemannian manifold the convex hull $\mathfrak{C}$ of $Q$ is strictly larger than  $Q_i$, for all $i$.

\parit{Proof of \ref{cor:caratheodory}.}
Without loss of generality we can assume that $\mathfrak{C}$ is a proper subset of $M$; in particular $\partial\mathfrak{C}\ne \varnothing$.
Since $Q$ is not contained in a geodesic, by the main theorem, $\mathfrak{C}$ has a non-trivial interior.
By the construction of $\mathfrak{C}$ above, any point $x\in \mathfrak{C} \setminus Q$ is not an extreme point of~$\mathfrak{C}$.

Assume $\partial \mathfrak{C} \not\subset Q$.
By the main theorem, the topological manifold $\partial \mathfrak{C}$ is the union of the closed subset $Q\cap \partial \mathfrak{C}$ and a countable union of geodesics.
But $\partial \mathfrak{C} \setminus Q$ is an $(m-1)$-dimensional topological manifold.
By dimensional reasons,  it is not a union of countably many rectifiable  curves --- a contradiction.
\qeds

\appendix

\section{Normalization of metrics}
\label{sec:normalization}

This appendix is devoted to the algebra of curvature tensor and its covariant derivatives that leads to a proof of Claim \ref{clm:codim-sigma}.

Choose an $m$-dimensional Euclidean space $\T$.
Denote by $\mathcal{S}$ the space of self-adjoint operators on $\T$.

Consider the space $\mathcal{G}$ of germs of Riemannian metrics on $\T$ at $0$ that coincide with the canonical metric at $0$.
Any germ in $\mathcal{G}$ can be described by $\langle G\cdot \vec v,\vec w\rangle$, where $\vec x\mapsto G_\vec x$ is a smooth function $\T\to\mathcal{S}$ such that $G_0=\id$. 

The $k$-jet of $G$ is defined by the Taylor polynomial of $G$ of degree $k$
$$G_\vec x=\id + G^1_\vec x+\dots+G^k_\vec x + o( |\vec x|^k ),
\eqlbl{eq:k-jet}$$
where $\vec x\mapsto G^i_\vec x$ is a homogeneous polynomial $\T\to\mathcal{S}$ of degree $i$.

Note that every array of homogeneous polynomials $G^1,\dots, G^k\:\T\to \mathcal{S}$ such that $\deg G^i\z=i$ appears in \ref{eq:k-jet} for the germ in $\mathcal{G}$ defined by 
\[G_\vec x=\id + G^1_\vec x+\dots+G^k_\vec x.
\eqlbl{eq:k-jet+}\]
The space of $k$-jets of germs in $\mathcal{G}$ will be denoted by~$\mathcal{G}^k$.

A germ in $\mathcal{G}$ will be called \emph{normal} if the standard coordinates on $\T$ coincide with normal coordinates of the germ in a neighborhood of the origin.
By the Gauss lemma, a germ defined by $G$ is normal if and only if 
\[G_\vec x\cdot \vec x=\vec x\eqlbl{eq:Gauss}\]
for all small $\vec x\in \T$.
The subspace of normal germs in $\mathcal{G}$ and their $k$-jets will be denoted by $\mathcal{N}$ and $\mathcal{N}^k$, respectively.

Suppose that $G$ describes a germ in $\mathcal{N}$
and $G^1,\dots, G^k$ be as in \ref{eq:k-jet}.
By \ref{eq:Gauss} 
\[G^i_\vec x\cdot \vec x=0
\eqlbl{eq:Gaussi}\] 
for any $i$.
Moreover, for an array of polynomials $G^1,\dots,G^k\:\T\to\mathcal{S}$ such that $G^i$ is homogeneous of degree $i$ and \ref{eq:Gaussi} holds for each $i$, the sum \ref{eq:k-jet+} defines a normal $k$-jet;
that is, \ref{eq:Gaussi} is the only condition on the normality of jets.

Christoffel symbols vanish in normal coordinates, thus, $G^1=0$, for  $G\in \mathcal N$.

Choose $\vec x\in \T$; denote by $\mathcal{S}_\vec x$ the subspace of the operators $S\in\mathcal{S}$ such that $S\cdot \vec x=0$.
By \ref{eq:Gaussi}, $G^i_\vec x\in \mathcal{S}_\vec x$ for any germ in $\mathcal{N}$.
The following claim says that $G^i_\vec x$ can be chosen arbitrarily in $\mathcal{S}_\vec x$ for $i\ge 2$ and $\vec x\ne 0$.

\begin{thm}{Claim}\label{clm:allSX}
Given $\vec x\ne 0$ in $\T$ and a sequence of operators $A_2,\dots A_k\in \mathcal{S}_\vec x$ there is a germ $(G^1,\dots,G^k)$ in $\mathcal{N}^k$ such that $G^i_\vec x=A_i$ for any $i\ge 2$.
\end{thm}

\parit{Proof.}
For any unit vector  $\vec y$ in $\T$  perpendicular to $\vec x$, consider the orthogonal projection $P^{\vec y}$
in $T$ onto the line generated by $\vec y$.  Diagonalizing operators  in $\mathcal S _{\vec x}$, we see that
such projections $P^{\vec y}$ generate $\mathcal S _{\vec x}$ as a vector space. 

$\mathcal N^k$ is described by \ref{eq:Gaussi}, hence it defines a linear subspace of $\mathcal S _{\vec x} ^k$.
Thus, it   suffices to  verify the following: For any $2\leq j \leq k$ and 
any unit vector  $\vec y$ in $\T$  perpendicular to $\vec x$, there exists a germ $(G^1,\dots,G^k)$ in $\mathcal{N}^k$ such that $G^j_\vec x=P^{\vec y}$ and $G^i =0$ for $i\neq j$.

 Such a normal germ can be constructed as a product 
of a surface of revolution (corresponding to the $(\vec x ,\vec y)$-plane) and a Euclidean space.
%
\qeds

Suppose that a germ in $\mathcal{G}$ is described by $G\:\T\to\mathcal{S}$.
Consider its array of Jacobi operators $(R^1,\dots,R^k)$ at the origin;
recall that $R^1=0$.
The identities in Section~\ref{sec:key} imply that any such array $(R^1,\dots,R^k)$ belongs to the space $\mathcal{R}^k$ defined by the following conditions: (i)
each $R^i\:\T\to \mathcal{S}$ is a homogeneous polynomial,
(ii) $\deg R^i=i$,
and (iii) $R^i_\vec x\cdot \vec x=0$ for any $i$ and $\vec x\in\T$.
Note that these conditions are exactly the same as for $G^i$ in $\mathcal{N}^k$.
Therefore, $\mathcal{R}^k$ can be identified with $\mathcal{N}^k$, but we will keep separate notations for them.

The expression of the curvature tensor in terms of the metric and its derivatives defines a natural algebraic map 
$$\rho_k\:\mathcal{G}^k\to \mathcal{R}^k.$$

For any $k\ge 2$, any $G\in \mathcal N^k$ and  $(R^1,\dots,R^k) = \rho_k (G)$,
we have
\[G^k=a_k\cdot  R^k + A^k,\eqlbl{eq:a+A}\]
where $a_k$ is a nonzero constant and $A^k$ is a field of self-adjoint operators that can be written as a polynomial of $R^2,\dots R^{k-2}$.
This statement follows easily from the formula derived by Old\v{r}ich Kowalski and Martin Belger \cite[Proposition 2.2]{kowalski-belger}.
(In fact, $a_k=-2\cdot\tfrac{k-1}{k+1}$, 
but we will not need it.)

Hence, the map $\rho _k$ admits an algebraic inverse map::

\begin{thm}{Claim}\label{clm:diff}
The restriction $\rho_k|_{\mathcal{N}^k}$ is an algebraic diffeomorphism $\mathcal{N}^k\leftrightarrow\mathcal{R}^k$.
\end{thm}

Applying \ref{clm:allSX}, we get the following:

\begin{thm}{Corollary}\label{cor:Rall}
Given $\vec x\ne 0$ in $\T$ and a sequence of operators $A_2,\dots A_k\in \mathcal{S}_\vec x$ there is a germ $\mathcal{N}^k$ with Jacobi operators $R^i_\vec x=A_i$ for any $i\ge 2$.
\end{thm}

\begin{thm}{Proposition}\label{prop:submersion}
The map $\rho_k:\mathcal{G}^k\to \mathcal{R}^k$ is an algebraic submersion.
\end{thm}

\begin{wrapfigure}{r}{35mm}
\vskip-4mm
\centering
\begin{tikzpicture}[scale=2]
\node (1) at (0,1) {$\mathcal{G}$};
\node (2) at (1.4,1){$\mathcal{N}$};
\node (11) at (0,0){$\mathcal{G}^k$};
\node (12) at (1.4,0) {$\mathcal{N}^k$};
\node (21) at (.7,-5/6) {$\mathcal{R}^k$};
\draw[>=latex, auto=right, loop above/.style={out=75,in=105,loop}, every loop,]
   (1) edge[bend left] node [swap] {$\nu$}  (2)
   (2) edge[bend left, left hook-stealth] node {$\iota$} (1)
   (1) edge (11)
   (2) edge (12)
   (11) edge[bend right] node  {$\rho_k$} (21)
   (11) edge [bend left] node [swap] {$\nu_k$} (12)
   (12) edge[bend left, left hook-stealth] node {$\iota_k$} (11)
    (21) edge[<->, bend right] node  {$\rho_k$} (12)
   ;
\end{tikzpicture}
\label{diagram-page}
\end{wrapfigure}

\parit{Proof.}
Evidently, $\rho_k$ is algebraic.

Any germ in $\mathcal{G}$ becomes normal if the space $\T$ is reparametrized by its exponential map.
This defines the normalization map 
$\mathcal{G}\xrightarrow{\nu} \mathcal{N}$.   Since the curvature tensors does not change under this (or any other) coordinate
change, $\nu$ commutes with   $\rho_k : \mathcal {G}, \mathcal {N} \to \mathcal R^k$.


By \ref{clm:diff}, $\mathcal{N}^k\xleftrightarrow{\rho_k}\mathcal{R}^k$ is a diffeomorphism.
The maps $\mathcal{G}^k\xrightarrow{\rho_k}\mathcal{R}^k\z{\xleftrightarrow{\rho_k}}\mathcal{N}^k$ together with the forgetful maps $\mathcal{G}\to\mathcal{G}^k$ and $\mathcal{N}\to\mathcal{N}^k$  commute.
In particular, we get a map $\mathcal{G}^k\xrightarrow{\nu_k}\mathcal{N}^k$ that commutes with the forgetful maps and the normalization $\nu$.  Hence, $\nu _k$ and $\rho _k$ commute.

Note that the inclusion $\mathcal{N}\stackrel{\iota}{\hookrightarrow} \mathcal{G}$ is a right inverse of~$\nu$.
Moreover, by changing the parametrization on $\T$ to normal coordinates of a given germ $G$ in $\mathcal{G}$, we may assume that $G$ lies in the image of~$\iota$.
Therefore, there is an inclusion $\mathcal{N}^k\stackrel{\iota_k}{\hookrightarrow} \mathcal{G}^k$ that is a right inverse of $\nu$ such that its image contains any given jet  in $\mathcal{G}^k$.
It follows that $\mathcal{G}^k\xrightarrow{\nu^k} \mathcal{N}^k$ is a submersion, hence the result.
\qeds

\parit{Proof of \ref{clm:codim-sigma}.}
Denote by $\tilde{\mathcal{G}}^k$ the space of all $k$-jets of Riemannian metrics at a given point $p$.
Denote by $\tilde\Sigma^k$ all jets in $\tilde{\mathcal{G}}^k$ such that for some nonzero tangent vector $\vec x\in\T_p$ the Jacobi operators $R^2_\vec x,\dots,R^k_\vec x$ have a common exceptional invariant subspace.

By Tarski--Seidenberg theorem, $\tilde\Sigma^k$ is semialgebraic; in particular it is stratified.
Due to the Thom transversality theorem \cite[2.3.2]{eliashberg-mishachev}, it is sufficient to show that for any point $p$ the codimension of
$\tilde\Sigma^k$  in $\tilde{\mathcal G}^k$ is larger than $m=\dim M$.

This is a pointwise statement;
therefore we may fix $p$ from now on.

A jet in $\tilde{\mathcal{G}}^0$
is described by the metric tensor $g_0$  on $\T=\T_pM$.
Note that the forgetful map $\tilde{\mathcal{G}} ^k \to \tilde{\mathcal{G}}^0$ is a fiber bundle.
Furthermore, the restriction of this forgetful map to $\tilde \Sigma$ is also a fiber bundle.
Thus, it suffices to prove that the intersection $\Sigma^k$ of $\tilde \Sigma^k$ with a fiber of $\tilde{\mathcal{G}} ^k \to \tilde{\mathcal{G}}^0$
has codimension at least $m$.
Note that the fiber of the forgetful map over the  Euclidean structure
on $\T$ given by $g_0$ is exactly the space $\mathcal G^k$ investigated above.

In other words, if we choose a chart $\T\to M$, then $g_0$ defines an inclusion $\mathcal{G}^k\hookrightarrow\tilde{\mathcal{G}}^k$,
and it is sufficient to show that 
\[\codim\Sigma^k\to\infty
\quad\text{as}\quad
k\to \infty;
\eqlbl{eq:codim-Sigma}\]
here we consider $\Sigma^k\z=\tilde\Sigma^k\cap \mathcal{G}^k$ as a subset of $\mathcal{G}^k$.

Denote by $\mathcal{L}$ the semialgebraic set of all pairs $(L,\vec x)$ where $L$ is a subspace of $\T$ such that $1<\dim L<m$,  and $\vec x\in L\backslash \{0\}$.
Given $(L,\vec x)\in\mathcal{L}$, denote by $\Sigma^k(L,\vec x)$ the subset of jets in $\mathcal{G}^k$ such that $L$ is an invariant subspace of all Jacobi operators $R^i_\vec x$  for any $i\le k$.

Choose $(L,\vec x)\in\mathcal{L}$. We claim that
\ref{prop:submersion} implies
\[\codim\Sigma^k(L,\vec x)\to\infty\quad \text{as}\quad k\to \infty.\eqlbl{eq:codim-Sigma(L,x)}\]
Indeed, a normal germ $(G^1,\dots G^k)$ belongs to $\Sigma^k(L,\vec x)$ if and only if all the Jacobi operators 
$R^2_\vec x,\dots, R^k_\vec x\in \mathcal{S}_\vec x$ have invariant subspace $L$.
The codimension of 
the space of $(k-1)$-tuples in $S_x$ that all have $L$ as invariant subspace grows with~$k$.
By \ref{prop:submersion} and \ref{cor:Rall} the composition 
$\mathcal{G}^k\z\to \mathcal{R}^k\z\to\mathcal{S}^{k-1}_\vec x$ that sends a germ to the array of its Jacobi operators
$(R^2_\vec x,\dots, R^k_\vec x)$ is a submersion.
Therefore, \ref{eq:codim-Sigma(L,x)} follows.

Observe that 
\[\codim\Sigma^k\ge \codim\Sigma^k(L,\vec x)-\dim\mathcal{L}.\]
Therefore, \ref{eq:codim-Sigma} follows.
\qeds

\section{Final remarks}\label{app:remarks}

We expect that the following question admits an affirmative answer.

\begin{thm}{Question}
Is it true that any Riemannian manifold $(M,g)$ contains a nontrivial geodesic  that runs in the boundary of some convex 
subset?
\end{thm}

There is a good chance that the argument of Albert Borb\'{e}ly  \cite[Lemma 2.1]{borbely} can be modified to answer the following question. 
Assuming that the answer is affirmative, it can be combined with the main proposition to derive further restrictions on convex hulls in generic Riemannian manifolds.

\begin{thm}{Question}
Let $\mathfrak{C}$ be the closure of a convex hull of a set $Q$ in a Riemannian manifold.
Then all points of $\mathfrak {C}$ with rank at most $1$ lie on minimizing geodesics between points in $Q$.
\end{thm}

The presented argument, when properly extended to infinite-dimensional manifolds, might lead to a negative answer to the following question of Mikhael Gromov \cite[6.B$_1\text{(f)}$]{gromov-1993}.

\begin{thm}{Question}
Let $X$ be a complete  $\CAT(0)$ space (not necessarily locally compact).
Is it true that any compact set of $X$ lies in a compact convex subset?
\end{thm}

A surprising behavior of convex sets in complete (but not locally compact) $\CAT(0)$ spaces is discussed by Nicolas Monod \cite{monod}.

Finally let us mention that there is a result of Anatoliy Milka \cite[§~4]{milka} about rank of points on geodesics in the intrinsic metric of convex surfaces; it is closely related to our main proposition but goes in the opposite direction.

{\sloppy
\printbibliography[heading=bibintoc]
\fussy
}
\end{document}